\documentclass[11pt,a4paper]{amsart}
\usepackage[left=3.5cm,right=3.5cm,top=3.5cm,bottom=3.5cm]{geometry}
\usepackage{enumerate}
\usepackage{leftidx}
\usepackage{mathtools}
\usepackage{amsmath,amscd,amssymb,amsthm}
\usepackage[arrow,matrix]{xy}
\usepackage[OT2,T1]{fontenc}
\usepackage{tikz}
\usepackage{graphicx}
\usepackage{hyperref}

\DeclareMathOperator{\im}{Im}
\DeclareMathOperator{\aut}{Aut}
\DeclareMathOperator{\ab}{ab}
\DeclareMathOperator{\GL}{GL}

\theoremstyle{definition}
\newtheorem{definition}{Definition}[section]

\newtheorem*{remark*}{Remark}
\newtheorem*{theorem*}{Theorem}

\theoremstyle{plain}
\newtheorem{theorem}[definition]{Theorem}
\newtheorem{corollary}[definition]{Corollary}
\newtheorem{lemma}[definition]{Lemma}
\newtheorem{proposition}[definition]{Proposition}
\newtheorem{conjecture}[definition]{Conjecture}

\newcommand{\F}{{ \mathbb F }}
\newcommand{\N}{{ \mathbb N }}
\newcommand{\Q}{{ \mathbb Q }}
\newcommand{\Z}{{ \mathbb Z }}
\newcommand{\R}{{ \mathbb R }}
\def \C {{\mathbb C}}

\renewcommand{\P}{{ \mathbb P }}

\def\house#1{{%
    \setbox0=\hbox{$#1$}
    \vrule height \dimexpr\ht0+1.4pt width .5pt depth \dp0\relax
    \vrule height \dimexpr\ht0+1.4pt width \dimexpr\wd0+2pt depth \dimexpr-\ht0-1pt\relax
    \llap{$#1$\kern1pt}
    \vrule height \dimexpr\ht0+1.4pt width .5pt depth \dp0\relax}}

\author[A. Ferraguti]{Andrea Ferraguti}
\address{Scuola Normale Superiore, Piazza dei Cavalieri 7, 56126 Pisa}
\email{and.ferraguti@gmail.com}

\title{A survey on abelian dynamical Galois groups}

\begin{document}

\begin{abstract}
Let $K$ be a number field, $f\in K[x]$ and $\alpha\in K$. A recent conjecture of Andrews and Petsche predicts that the dynamical Galois group of the pair $(f,\alpha)$ is abelian if and only if the pair $(f,\alpha)$ is $K^{\ab}$-conjugated to $(g,\beta)$, where $g$ is a power or a Chebyshev map and $\beta$ is $\zeta$ or $\zeta+\zeta^{-1}$, respectively, and $\zeta$ is a root of unity. We review three completely different approaches that allow to prove several cases of the conjecture.

\end{abstract}

\maketitle

\section{Introduction}

Let $K$ be a field, $f\in K(x)$ be a rational map of degree at least $2$ and $\alpha\in K$. For every $n\geq 1$ we will denote by $f^n$ the $n$-th iterate of the map $f$, by $K_n(f,\alpha)$ the splitting field of $f^n-\alpha$, by $G_n(f,\alpha)$ its Galois group, by $K_\infty(f,\alpha)$ the direct limit of the $K_n(f,\alpha)$'s and by $G_\infty(f,\alpha)$ the inverse limit of the $G_n(f,\alpha)$ (or, equivalently, the Galois group of $K_\infty(f,\alpha)$). The group $G_\infty(f,\alpha)$ is called \emph{dynamical Galois group} of the pair $(f,\alpha)$.

The study of dynamical Galois groups is partially motivated by the fact that certain instances of them arise as central objects in the classical context of elliptic curves. Namely, if $E$ is an elliptic curve over a field $K$, $G_K$ is the absolute Galois group of $K$, and $p$ is a prime different from $\text{char } K$, the Galois representation $\rho_{E,p}$ attached to the $p$-adic Tate module of $E$ is a continuous homomorphism $G_K\to \aut(T_p(E))\cong \GL_2(\Z_p)$ whose image is essentially the dynamical Galois group of the pair $([p],\mathcal O_E)$, where $[p]$ is the multiplication-by-$p$ map and $\mathcal O_E$ is the identity for the group structure on $E$. It is interesting to notice that although the construction of dynamical Galois groups arising from self-rational maps of $\mathbb P^1$ can be generalized in a straightforward way to self-morphisms of any algebraic curve, there are no non-trivial examples where the curve is smooth and has genus at least $2$, by Riemann-Hurwitz formula. Hence $\mathbb P^1$ and elliptic curves are the only interesting smooth curves for studying these objects.

When the base field $K$ is a number field, one has the following celebrated result, due to Serre \cite{serre}.
\begin{theorem*}
The index $[\GL_2(\Z_p):\im(\rho_{E,p})]$ is infinite if and only if $E$ has complex multiplication.
\end{theorem*}

This is part of a general philosophy according to which whenever we have a continuous action of $G_K$ on a topological vector space $V$ (or, equivalently, a continuous representation $G_K\to \aut(V)$) the existence of "extra symmetries" on $V$ has the effect of constraining the image of the associated representation to be "small". The introduction of \cite{conti} contains a nice account of this idea for Galois representations associated to various arithmetic objects.

When dealing with dynamical Galois groups arising from self-rational maps of $\mathbb P^1$ one expects a similar beaviour. If $f\in K(x)$ and $\alpha\in K$ is such that for every positive integer $n$ the set $f^{-n}(\alpha)\coloneqq (f^n)^{-1}(\alpha)$ contains no ramification points for $f$, then $G_\infty(f,\alpha)$ can be seen as the image of the arboreal representation associated to the pair $(f,\alpha)$. This is a continuous homomorphism $\rho_{f,\alpha} \colon G_K\to \aut(T_\infty(f,\alpha))$, where $T_\infty(f,\alpha)$ is the infinite regular tree constructed with the \emph{backward orbit} of $\alpha$, i.e.\ the set $f^{-\infty}(\alpha)\coloneqq \bigcup_{n\ge 1}f^{-n}(\alpha)$ (see \cite{jones_survey} for details). This framework is far poorer of structure than the elliptic curve one; here for example the map $\rho_{f,\alpha}$ does not arise from the action of $G_K$ on a vector space but it is just a continuous homomorphism of topological groups. Nevertheless, the pair $(f,\alpha)$ can admit extra symmetries that force the image of $\rho_{f,\alpha}$ to be small, i.e.\ to have infinite index image in $\aut(T_\infty(f,\alpha))$.

When $K$ is a global field of characteristic not $2$ and $f$ is a quadratic map, a conjecture of Jones \cite[Conjecture 3.11]{jones_survey} describes all possible extra symmetries that force the index $[\aut(T_\infty(f,\alpha)):\im(\rho_{f,\alpha})]$ to be infinite. When restricted to quadratic polynomials, the conjecture is the following. Recall that a rational map $f\in K(x)$ is called \emph{post-critically finite} (PCF for short) if the orbits of its ramification points, called \emph{critical points} by dynamicists, are all finite.
\begin{conjecture}\label{jones_conjecture}
Let $K$ be a global field of characteristic not $2$, let $f\in K[x]$ be a quadratic polynomial and $\alpha\in K$. Suppose that $f^{-\infty}(\alpha)$ contains no ramification points for $f$. Then the index $[\aut(T_\infty(f,\alpha)):\im(\rho_{f,\alpha})]$ is infinite if and only if $f$ is post-critically finite or $\alpha$ is periodic for $f$.
\end{conjecture}
 The "if" part of the conjecture is actually a theorem (see \cite{jones_survey}); the difficult part is showing necessity. Some results conditional on the $abc$ conjecture for number fields are proven in \cite{juul}.

Notice that the two types of symmetries appearing in Conjecture \ref{jones_conjecture} are quite different in nature: being post-critically finite is a property of the polynomial $f$: no matter how we choose the basepoint $\alpha$ the group $G_\infty(f,\alpha)$ is going to be small. This means that the so-called \emph{iterated monodromy group} of the map $f$, namely the dynamical Galois group of the pair $(f,t)$ where $t$ is transcendental over $K$, is small as well. On the other hand, $\alpha$ being periodic for $f$ is a symmetry of the pair $(f,\alpha)$. The iterated monodromy group of $f$ might even be large, for example if $f=x^2-3$ and $K=\Q$ we have that $G_\infty(f,t)\cong \aut(T_\infty(f,\alpha))$ (see \cite{pink} for the general theory of iterated monodromy groups of quadratic functions), but specializing at $t=1$ makes the dynamical Galois group of the specialized pair far smaller. One would expect that most specializations do not have such a drastic impact on the size of the dynamical Galois group. Certainly $f$ has only a finite number of $K$-rational periodic points, for example. In \cite{benedetto1}, the authors propose a conjecture along these lines for the iterated monodromy group of $x^2-1$.

For all the reasons above, it seems natural to expect that dynamical Galois groups can rarely be abelian. In the elliptic curves framework, abelian images of Galois representations are fully described by complex multiplication theory. Let $K$ be a number field, $E/K$ be an elliptic curve, $p$ a prime and $\rho_{E,p}$ the Galois representation of $G_K$ on the $p$-adic Tate module of $E$. Then the following theorem holds (see for example \cite{lozano}).
\begin{theorem*}
The image of $\rho_{E,p}$ is virtually abelian if and only if $E$ has complex multiplication, and it is abelian if and only if $E$ has complex multiplication defined over $K$.
\end{theorem*}

That is, in the elliptic curve case the existence of extra symmetries not only has the effect of forcing the image of the associated representation to be small, it even forces it to be (virtually) abelian.

Once again, the dynamical side of the story is far more mysterious, and only recently a conjecture has been formulated in \cite{andrews}. Recall that if $f,g\in K(x)$, $\alpha,\beta\in K$ and $L$ is an extension of $K$, we say that $(f,\alpha)$ and $(g,\beta)$ are \emph{conjugated} over $L$ if there exist a M\"obius transformation $m\in \text{PGL}_2(L)$ such that $g=m^{-1}\circ f\circ m$ and $m(\beta)=\alpha$. If the two pairs $(f,\alpha)$ and $(g,\beta)$ are $K$-conjugated, the two groups $G_\infty(f,\alpha)$ and $G_\infty(g,\beta)$ are isomorphic. Moreover, a point $\alpha\in \overline{K}$ is \emph{exceptional} for $f$ if its backward orbit is finite.

\begin{conjecture}[{{\cite[Conjecture]{andrews}}}]\label{ap_conjecture}
Let $K$ be a number field, $f\in K[x]$ a polynomial of degree $d\ge 2$ and $\alpha \in K$ be a non-exceptional point for $f$. Then $G_\infty(f,\alpha)$ is abelian if and only if there exists a root of unity $\zeta$ such that the pair $(f,\alpha)$ is $K^{\text{ab}}$-conjugated to either $(x^d,\zeta)$ or $(\pm T_d(x),\zeta+\zeta^{-1})$, where $T_d$ is the Chebyshev polynomial of degree $d$.\footnote{The formulation of \cite{andrews} overlooked the case of $-T_d(x)$. Notice that $d$ and $T_d$ have the same parity, and hence it is immediate to see that $G_\infty(-T_d(x),\zeta+\zeta^{-1})$ is abelian. However, $-T_d(x)$ and $T_d(x)$ are never conjugated over $\overline{K}$ when $d$ is odd.}
\end{conjecture}

Hence in the dynamical framework the current prediction is that abelian dynamical Galois groups can only occur in a very special case: not only the polynomial needs to be PCF, but it also needs to come from an endomorphism of an algebraic group.

The sufficiency part of Conjecture \ref{ap_conjecture} has been proven by the same authors, as below.

\begin{theorem}[{{\cite[Theorems 12 and 13]{andrews}}}]\label{ap_sufficiency}
Let $K$ be a number field, $f\in K[x]$ be a polynomial of degree $d\ge 2$ and $\alpha\in K$ be a non-exceptional point for $f$. Suppose that $f$ is $\overline{K}$-conjugate to $x^d$ or $\pm T_d(x)$. Then $G_\infty(f,\alpha)$ is abelian if and only if $(f,\alpha)$ is $K^{\text{ab}}$-conjugated to $(x^d,\zeta)$ or $(\pm T_d(x),\zeta+\zeta^{-1})$, for some root of unity $\zeta$.
\end{theorem}

This survey is devoted to illustrate three approaches, based on very different strategies, that have been used to prove several cases of Conjecture \ref{ap_conjecture}. Notice that since power and Chebyshev maps are both PCF, a natural idea in order to prove Conjecture \ref{ap_conjecture} is to start by proving that if $G_\infty(f,\alpha)$ is abelian then $f$ is PCF. If the base field $K$ is fixed there are only finitely many $\overline{K}$-conjugacy classes of PCF rational maps of degree $d$, outside the flexible Latt\`es locus (see \cite{benedetto2}), that contains no polynomials. Hence showing the above implication already reduces considerably the complexity of the problem: for example when $K=\Q$ and $d=2$ the only PCF polynomials are (up to conjugacy over $\mathbb Q$) $x^2,x^2-1$ and $x^2-2$. The approaches described in sections \ref{height} and \ref{groups} both follow this route, while the third one uses a completely different strategy.

\section{Height pairing approach}\label{height}
In this section we illustrate the route taken by the authors of \cite{andrews}, which proved the following theorem. Recall that if $K$ is a field, $f\in K[x]$ and $\alpha\in K$ the pair $(f,\alpha)$ is said to be \emph{stable} if for every $n\geq 1$ the polynomial $f^n-\alpha$ is irreducible over $K$.

\begin{theorem}\label{ap_theorem}
Conjecture \ref{ap_conjecture} holds true when $K=\Q$, $f\in \Q[x]$ is a quadratic polynomial and $\alpha\in \Q$ is such that $(f,\alpha)$ is stable.
\end{theorem}
As explained in the Introduction, the first step is to be able to reduce to the PCF case. This is done in \cite[Theorem 8]{andrews}, where the authors prove that if $f\in K[x]$ is quadratic, $K$ is a number field, $(f,\alpha)$ is stable and $G_\infty(f,\alpha)$ is abelian, then $f$ is PCF. However, as we will see in Section \ref{groups}, it is possible to remove the stability hypothesis by following a different approach. Theorem \ref{ap_theorem} is superseeded by Theorems \ref{fp_theorem} and \ref{o_theorem}, but it is nevertheless interesting to understand the proof strategy, as it is completely different from that of the aforementioned theorems.

The PCF reduction argument implies that in order to prove Theorem \ref{ap_theorem} it is enough to show that if $f\in \{x^2,x^2-1,x^2-2\}$ and $\alpha\in \Q$ is such that $(f,\alpha)$ is stable and $G_\infty(f,\alpha)$ is abelian, then Conjecture \ref{ap_conjecture} holds true. Notice that $x^2-2$ is the Chebyshev polynomial of degree 2. The hard part of the proof is to show that $G_\infty(x^2-1,\alpha)$ cannot be abelian.

In order to explain the proof strategy, let us start with the toy case $f=x^d\in K[x]$, with $d\ge 2$, $K$ a number field and $\alpha\in K^{\times}$ (notice that $0$ is an exceptional point for $x^d$). Theorem \ref{ap_sufficiency} says that $G_\infty(x^d,\alpha)$ is abelian if and only if $\alpha$ is a root of unity. Sufficiency is obvious, so now assume that $G_\infty(x^d,\alpha)$ is abelian. This is equivalent to saying that the whole backward orbit of $\alpha$ is contained in $K^{\text{ab}}$. Now observe that if $\gamma_n\in \overline{K}$ is such that $f^n(\gamma_n)=\alpha$ for some $n\geq 1$, then
$$h(\gamma_n)=\frac{1}{d^n}h(\alpha),$$
where $h(\cdot)$ is the absolute logarithmic Weil height. Hence if $\{\gamma_n\}_{n\geq 0}\subseteq \overline{K}$ is a sequence of algebraic numbers with $\gamma_0=\alpha$ and $f^n(\gamma_n)=\alpha$ for every $n\ge 1$, then $h(\gamma_n)\to 0$ as $n\to +\infty$. Now the key tool is the following theorem of Amoroso and Zannier \cite{amoroso}.
\begin{theorem}\label{height_bound}
Let $K$ be a number field. There exists a positive constant $c_K$ such that if $\gamma\in K^{\text{ab}}$ is not zero nor a root of unity, then $h(\gamma)\geq c_K$.
\end{theorem}
Hence the fact that all the $\gamma_n$'s belong to $K^{\ab}$ and their heights tend to zero implies that they must all be roots of unity, and in turn $\alpha$ must be a root of unity.

Of course the argument for $\pm T_d(x)$ is also very similar, using the fact that $T_d(x+x^{-1})=x^d+x^{-d}$.

Now the idea for proving that if $\alpha\in \Q$ and $(x^2-1,\alpha)$ is stable then $G_\infty(x^2-1,\alpha)$ cannot be abelian is similar in spirit, but it uses more sophisticated tools. The first one is the so-called \emph{Call-Silverman canonical height} (see \cite[Chapter 3]{silverman2}), which is a function $\widehat{h}_\phi\colon\overline{K}\to \R_{\geq 0}$ associated to a polynomial $\phi\in K[x]$ where $K$ is a number field. This is defined as follows: let $\phi\in K[x]$ be a degree $d$ polynomial and let $\gamma\in \overline{K}$. Then
$$\widehat{h}_\phi(\gamma)\coloneqq\lim_{n\to \infty}\frac{h(\phi^n(\gamma))}{d^n}.$$
This can be viewed as a dynamical generalization of the absolute logarithmic Weil height, and in fact it differs from it by a bounded amount. Notice that on the one hand we have that $\widehat{h}_{x^d}=h$, and on the other one $\widehat{h}_{\phi}(\phi(\gamma))=d\cdot \widehat{h}_\phi(\gamma)$. An algebraic number $\gamma$ satisfies $\widehat{h}_\phi(\gamma)=0$ if and only if it is preperiodic for $\phi$; this is consistent with the fact that preperiodic points of $x^d$ are exactly $0$ and roots of unity. A crucial property of the Call-Silverman height is that it decomposes as a (weighted) sum of local components $\lambda_{\phi,v}$ for each place $v$ of $K$.

The next crucial ingredient in the proof is the \emph{Arakelov-Zhang pairing}. This is defined, for two polynomials $\phi,\psi\in K[x]$ as
$$\langle\phi,\psi\rangle\coloneqq \sum_{v}\int_{\C_v}\lambda_{\psi,v}d\mu_{\phi,v}.$$
Here the sum is over all places $v$ of $K$, while $\mu_{\phi,v}$ is the canonical measure associated to $\phi$, and is a $\phi$-invariant Borel unit measure supported on the filled Julia set of $\phi$. It is not important for the purpose of this paper to go into detail about it, the interested reader can refer to \cite{baker,chambert,petsche}. The Arakelov-Zhang pairing has three key properties: it is non-negative, it is symmetric and if $\{\gamma_n\}_{n\in \N}\subseteq \overline{K}$ is a sequence of points such that $\widehat{h}_\phi(\gamma_n)\to 0$, then $\widehat{h}_\psi(\gamma_n)\to \langle \phi,\psi\rangle$ (see \cite{petsche}).

The properties of the Call-Silverman height and the Arakelov-Zhang pairing suggest a strategy to prove Theorem \ref{ap_theorem}. In fact, suppose that $\alpha \in \Q$ is such that $G_\infty(x^2-1,\alpha)$ is abelian, and let $\{\gamma_n\}_{n\geq 0}\subseteq \overline{\Q}$ be a sequence of points such that $\gamma_0=\alpha$ and $\gamma_n^2-1=\gamma_{n-1}$ for every $n\geq 1$. Then $\widehat{h}_{x^2-1}(\gamma_n)\to 0$, and therefore $h(\gamma_n)=\widehat{h}_{x^2}(\gamma_n)\to \langle x^2-1,x^2\rangle$, that is a positive constant (see \cite{petsche}) and it can be numerically approximated, as explained in \cite{andrews}. On the other hand, since all $\gamma_n$'s belong to $\Q^{\ab}$, their logarithmic Weil height must be either zero or at least the constant $c_\Q$ arising from Theorem \ref{height_bound}. If now we are lucky enough that $\langle x^2,x^2-1\rangle <c_\Q$, the theorem is proven. Unfortunately, this is not the case. However, the hypothesis under which Andrews-Petsche work, namely the fact that $(x^2-1,\alpha)$ is stable, allows to show, using ramification theory, that if $G_\infty(x^2-1,\alpha)$ is abelian then $\Q_\infty(x^2-1,\alpha)\subseteq \Q(\zeta_{2^\infty})$. Clearly since Theorem \ref{height_bound} holds for $\Q^{\ab}$ then it holds for $\Q(\zeta_{2^\infty})$ with a larger constant. It turns out that such constant is computable, and is indeed larger than $\langle x^2,x^2-1\rangle$, completing the proof of Theorem \ref{ap_theorem}.

\section{Group theoretic approach}\label{groups}

The approach taken by the author and Pagano in \cite{ferra1} is of a completely different nature than the one described in Section \ref{height}. In fact, it is based on the group-theoretical properties of $\aut(T_\infty(f,\alpha))$ when $f$ is quadratic. Let $\Omega_\infty$ denote, from now on, the group of automorphisms of the infinite, regular, rooted binary tree (see \cite{ferra2} for details about the structure of the group). The reason why the quadratic case differs heavily from the general one is that $\Omega_\infty$ is a pro-2-group, and therefore its structure is far more accessible than in the general case.

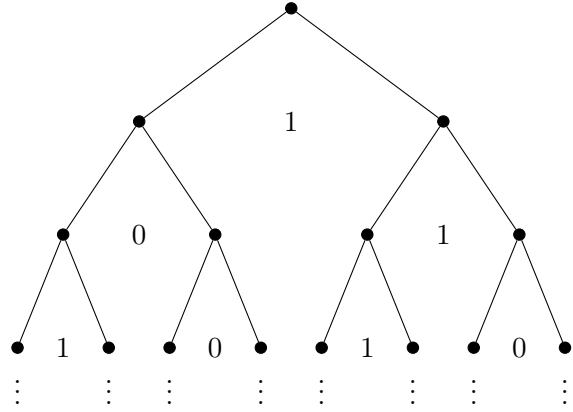
\begin{figure}
\centering

\begin{tikzpicture}
\tikzstyle{solid node}=[circle,draw,inner sep=1.5,fill=black]
\node[solid node]{}[sibling distance=40mm]
child{node[solid node](a){}[sibling distance=20mm]
child{node[solid node](c){}[sibling distance=12mm]
child{node[solid node,label=below:\vdots](g){}}
child{node[solid node,label=below:\vdots](h){}}}
child{node[solid node](d){}[sibling distance=12mm]
child{node[solid node,label=below:\vdots](i){}}
child{node[solid node,label=below:\vdots](l){}}}
}
child{node[solid node](b){}[sibling distance=20mm]
child{node[solid node](e){}[sibling distance=12mm]
child{node[solid node,label=below:\vdots](m){}}
child{node[solid node,label=below:\vdots](n){}}}
child{node[solid node](f){}[sibling distance=12mm]
child{node[solid node,label=below:\vdots](o){}}
child{node[solid node,label=below:\vdots](p){}}}
};
\path (a) -- (b) node [midway] {1};
\path (c) -- (d) node [midway] {0};
\path (e) -- (f) node [midway] {1};
\path (g) -- (h) node [midway] {1};
\path (i) -- (l) node [midway] {0};
\path (m) -- (n) node [midway] {1};
\path (o) -- (p) node [midway] {0};
\end{tikzpicture}\caption{The portrait of an element $\sigma=(1,01,1010,\ldots)$. If we label from $1$ to $8$, from left to right, the nodes at distance $3$ from the root, $\sigma$ corresponds to the permutation $(1736)(2845)$ on those nodes.}\label{fig1}
\end{figure}

Elements of $\Omega_\infty$ can be uniquely identified by their \emph{portrait}. If $\sigma\in \Omega_\infty$, its portrait is an infinite sequence of the following form:
$$(\varepsilon_1^1,\varepsilon_1^2\varepsilon_2^2,\varepsilon_1^3\varepsilon_2^3\varepsilon_3^3\varepsilon_4^3,\ldots,\varepsilon_1^{n}\ldots\varepsilon_{2^{n-1}}^{n},\ldots),$$
where $\varepsilon_i^j\in \F_2$ for every $i,j$. This has a very clear combinatorial interpretation: at each level $n$ of the tree, if the $2^{n-1}$ pairs of brother nodes are labeled from $1$ to $2^{n-1}$, the symbol $\varepsilon_i^n$ determines whether the nodes of the $i$-th pair are swapped by $\sigma$ or not. If $\sigma=1$ they are swapped, otherwise they are not. See Figure \ref{fig1} for how to picture the portrait of an element $\sigma=(1,01,1010,\ldots)$.

The portrait is a very useful tool for thinking about elements of $\Omega_\infty$. In \cite{ferra1} we showed that for every $n\geq 1$ there is a continuous surjective homomorphism
$$\phi_n\colon \Omega_\infty\to \F_2$$
$$\sigma=(\varepsilon_1^1,\varepsilon_1^{2}\varepsilon_2^2,\ldots)\mapsto \sum_{i=1}^{2^{n-1}}\varepsilon_i^n,$$
and the $\F_2$-vector space $\hom_{\text{cont}}(\Omega_\infty,\F_2)$ is generated by the $\phi_n$'s. This allows for a compact description of maximal subgroups of $\Omega_\infty$: they are precisely those of the form $\ker \sum_{i\in I}\phi_i$, where $I\subseteq \N$ is a finite subset. Moreover, the commutator subgroup $[\Omega_\infty,\Omega_\infty]$ coincides with $\bigcap_{n\geq 1}\ker\phi_n$, and there is an isomorphism $\Omega_\infty^{\ab}\stackrel{\cong}{\to} \prod_{n\geq 1}\F_2$. The composite map $\psi\colon\Omega_\infty\twoheadrightarrow \Omega_\infty^{\ab}\to \prod_{n\geq 1}\F_2$ is given by $\sigma\mapsto (\phi_1(\sigma),\phi_2(\sigma),\ldots,\phi_n(\sigma),\ldots)$. This fact can be used to prove the following very useful lemma.

\begin{proposition}[{{\cite[Proposition 7.3]{ferra5}}}]\label{propA}
Let $\sigma,\tau\in \Omega_\infty$ and suppose that $\phi_1(\sigma)=1$. If $\psi(\sigma)$ and $\psi(\tau)$ are linearly independent vectors in the $\F_2$-vector space $\prod_{n\geq 1}\F_2$, then $\sigma$ and $\tau$ do not commute.
\end{proposition}
It is an immediate consequence that if $G\leq \Omega_\infty$ is abelian and there is at least one $\sigma\in G$ with $\phi_1(\sigma)=1$, then $\psi(G)$ is an $\F_2$-vector space of dimension $1$.

Clearly this forces a very restrictive condition on the nature of abelian subgroups of $\Omega_\infty$; the problem is how to translate it on a condition on $(f,\alpha)$ in order for $G_\infty(f,\alpha)$ to be abelian. From now on, suppose that $K$ is a field of characteristic not $2$, let $f\in K[x]$ be a monic, quadratic polynomial and $\alpha\in K$ be an element such that $f^{-\infty}(\alpha)$ contains no ramification points for $f$, so that the image of the arboreal representation $\rho_{f,\alpha}\colon G_K\to \aut(T_\infty(f,\alpha))$ can be identified with a closed subgroup of $\Omega_\infty$. Up to conjugation over $K$, we can always assume that $f=x^2+c$ for some $c\in K$. Recall that the \emph{adjusted post-critical orbit} of $f$ is the sequence $\{c_n\}_{n\geq 1}\subseteq K$ given by $c_1=-c$, $c_n=f(c_{n-1})$ for every $n\geq 2$. We denote by $\{c_{n,\alpha}\}$ the sequence defined by $c_{1,\alpha}=c_1+\alpha$ and $c_{n,\alpha}=c_n-\alpha$ for every $n\geq 2$. Now the link between the group theoretical properties of $\Omega_\infty$ and the structure of $\rho_{f,\alpha}$ is provided by the following proposition.

\begin{proposition}{{\cite[Proposition 4.2]{ferra5}}}\label{propB}
Let $I\subseteq \N$ be a finite subset. The image of $\rho_{f,\alpha}$ is contained in $\ker\left(\sum_{i\in I}\phi_i\right)$ if and only if $\prod_{i\in I}c_i\in K^2$.
\end{proposition}
Propositions \ref{propA} and \ref{propB} together imply immediately the following theorem.

\begin{theorem}\label{1dimensionality}
If $\im \rho_{f,\alpha}$ is an abelian subgroup of $\Omega_\infty$ and $f-\alpha$ is irreducible over $K$, then the sub-$\F_2$-space generated by the sequence $\{c_{n,\alpha}\}$ inside the $\F_2$-vector space $K^{\times}/{K^{\times}}^2$ is $1$-dimensional.
\end{theorem}
When $K$ is a number field the $1$-dimensionality condition given by Theorem \ref{1dimensionality} is a heavy constraint on the arithmetic of $f$. In fact it is a standard expectation in arithmetic dynamics that if the orbit of $0$ is infinite, i.e.\ if $f$ is \emph{ post-critically infinite}, then all but finitely many elements of the sequence $\{c_{n,\alpha}\}$ admit a primitive prime divisor with odd valuation, and this is a theorem (with certain extra assumptions) under the abc conjecture (see \cite{tucker}). A standard application of Faltings' theorem on rational points on curves of high genus shows easily that when $K$ is a number field such $1$-dimensionality condition implies that $f$ must be PCF. A couple of observations show that for the sake of this implication it is possible to remove the conditions on $f^{-\infty}(\alpha)$ not to contain ramification points and on $f-\alpha$ to be irreducible; one therefore obtains the following corollary.
\begin{corollary}\label{pcf}
Let $K$ be a number field, let $\alpha\in K$ and let $f\in K[x]$ be a quadratic polynomial. If $G_\infty(f,\alpha)$ is abelian, then $f$ is PCF.
\end{corollary}

 However one should absolutely not overlook the conclusion on the $1$-dimensionality of the span of $\{c_{n,\alpha}\}$ given by Theorem \ref{1dimensionality} for two reasons: in the first place, it is useful in proving that $G_\infty(f,\alpha)$ is not abelian for a given pair $(f,\alpha)$, and in the second place it is extremely versatile because it holds over any number field, and not just over $\Q$.

Hence when $K=\Q$ and $f\in \Q[x]$ is quadratic we are again reduced to study the three polynomials $x^2,x^2-2$ and $x^2-1$, but as we already explained in Section \ref{height} the only hard case is $x^2-1$. In \cite{ferra1} we took once again a completely different route than the one used in \cite{andrews}. Our idea was that in order for the group  $G_\infty(x^2-1,\alpha)$ to be abelian, it needs to be locally abelian at every prime. In particular since $G_\infty(x^2-1,\alpha)$ is a pro-$2$-group, if $p$ is an odd prime then local class field theory implies that the extension $(\Q_p)_\infty(x^2-1,\alpha)$ must be finitely ramified. A general result of Anderson et al.\ \cite[Lemma 7.1, Theorem 7.4]{anderson} shows that if $p$ is odd then for $(\Q_p)_\infty(x^2-1,\alpha)$ to be finitely ramified one needs $v_p(\alpha)=v_p(\alpha+1)=0$, where $v_p$ is the usual $p$-adic valuation. But clearly there are only finitely many rational $\alpha$'s that satisfy this condition at every odd prime: they are $0,-1,1,-2,-1/2$. Now one can rule these values out by ad hoc arguments, or by explicitly computing the Galois groups of the first few iterates. Let us remark that this is an example of the usefulness of Theorem \ref{1dimensionality}: for example when $\alpha=-1/2$ one has $c_{1,\alpha}=1/2=c_{2,\alpha}$ while $c_{3,\alpha}=-1/2$. Now clearly the subspace of $\Q^{\times}/{\Q^{\times}}^2$ spanned by $1/2$ and $-1/2$ has dimension $2$, and therefore $G_\infty(x^2-1,-1/2)$ cannot be abelian (in fact, Theorem \ref{1dimensionality} can also be stated at a finite level, so that in this case the group $G_3(x^2-1,-1/2)$ is already non-abelian).

All in all, these arguments prove the following theorem.

\begin{theorem}\label{fp_theorem}
Conjecture \ref{ap_conjecture} holds true for every pair $(f,\alpha)$ with $f\in \Q[x]$ quadratic and $\alpha\in \Q$.
\end{theorem}  

In a forthcoming work of the author and Pagano, we will in fact use the same approach to show that the same result holds true for quadratic polynomials over any quadratic field.

We also remark that thanks to the great generality of Theorem \ref{1dimensionality}, one can also prove the following result.

\begin{theorem}[{{\cite[Theorem 3.13]{ferra1}}}]
Let $K$ be a global function field of odd characteristic $p$, let $f\in K[x]$ be quadratic and $\alpha\in K$. Then $G_\infty(f,\alpha)$ is abelian if and only if $(f,\alpha)$ is $K$-conjugated to a pair $(g,\beta)$ with $g\in \overline{\F}_p[x]$ and $\beta\in \overline{\F}_p$.
\end{theorem}
\section{Unlikely intersection approach}

The idea behind this approach had been introduced by Dvornicich and Zannier in \cite{DZ}, where it was used to study cyclotomic preperiodic points of rational functions. Later in \cite{ostafe} Ostafe used it to study cyclotomic points in backward orbits of rational functions. This technique is quite powerful, and as we will see it implies the validity of Conjecture \ref{ap_conjecture} for $K=\Q$. The downside however is that at the present time there is no hope of generalizing it to number fields other than $\Q$.

In order to state the result, let us introduce some notation. Let $K$ be a number field, $\alpha\in K$ and $\phi\in K(x)$ be a rational function of degree $d\ge 2$. Let $K^c$ be the cyclotomic closure of $K$ inside a fixed algebraic closure $\overline{K}$. We say that the pair $(\phi,\alpha)$ is \emph{special} if it is $K^c$-conjugate to one of the following pairs: $\{(x^{\pm d},\zeta),(\pm T_d(x),\zeta+\zeta^{-1})\}$ where $\zeta$ is a root of unity. We recall that we denote by $\phi^{-\infty}(\alpha)$ the backward orbit of $\alpha$, namely the set of all $\gamma\in \overline{K}$ such that $\phi^n(\gamma)=\alpha$ for some $n\ge 0$.

\begin{theorem}[{{\cite[pp.\ 1933--1935]{ostafe}}}]\label{o_theorem}
Let $\phi=f/g\in K(x)$, where $f,g\in K[x]$ are coprime polynomials. Suppose that $\deg f>\deg g+1$ and that $(\phi,\alpha)$ is not a special pair. Then $\phi^{-\infty}(\alpha)\cap K^c$ is finite.
\end{theorem}
In order to illustrate the proof of this theorem, we need three preliminary results. The first one relies crucially on the hypothesis $\deg f>\deg g+1$ that, from a dynamical point of view, says that $\infty$ is a superattracting fixed point for $\phi$. Recall that the \emph{house} of an algebraic number $\gamma$ is defined as $\house{\gamma}\coloneqq\max_{\sigma\colon \Q(\gamma)\to \C}|\sigma(\gamma)|$.

\begin{lemma}[{{\cite[Corollaries 2.6 and 2.8]{ostafe}}}]\label{absolute_bounds}
Let $\gamma\in \phi^{-\infty}(\alpha)$.
\begin{enumerate}
\item There exists a positive constant $H$, depending on $K,\phi,\alpha$, such that $\house{\gamma}\leq H$.
\item There exists a positive integer $D$, depending on $K,\phi,\alpha$, such that $D\cdot \gamma$ is an algebraic integer.
\end{enumerate}

\end{lemma}

Let us briefly explain the idea behind the above lemma. Since $\infty$ is superattracting, whenever $v$ is a place of $K$ and we consider $\phi$ as a rational map $\phi_v\colon\P^1_{K_v}\to \P^1_{K_v}$ then $\infty$ stays superattracting (this fails as soon as $\infty$ is not superattracting, as if its multiplier is nonzero then for all but finitely many places $\infty$ becomes an indifferent point). Now let $v$ be an archimedean place of $K$. Since $\infty$ is a superattracting point for $\phi_v\colon \P^1_{K_v}\to \P^1_{K_v}$ then there exists a neighborhood $U_v$ of $\infty$ such that $\phi_v(U_v)\subseteq U_v$ (see \cite{benedetto5}). Hence it is clear that if $U_v$ is chosen to be small enough so that $\alpha\notin U_v$, no point $\gamma \in U_v$ can be mapped to $\alpha$ by any iterate of $\phi$. Of course lying outside of $U_v$, since the latter is a neighborhood of $\infty$, translates into a bound on the absolute value of $\gamma$. Since there is a bound for every archimedean embedding, there is also a bound on the house. For non-archimedean places the situation is essentially the same. It is easy to see that there is a finite set $S$ of places of $K$ such that every point in the backward orbit of $\alpha$ is an $S$-integer: it is enough for $S$ to contain all places where $\alpha$ or some coefficient of $f$ or $g$ has negative valuation. For primes in $S$ the same argument of the archimedean case applies (see \cite[Proposition 4.3]{benedetto5}) and leads to the existence of a $D$ as in Lemma \ref{absolute_bounds}. The constants $H$ and $D$ can be made explicit.

The second preliminary result is the following theorem, that is a generalization proven by Dvornicich-Zannier of a result of Loxton \cite{loxton}.

\begin{theorem}[{{\cite[Theorem L]{DZ}}}]\label{loxton}
Let $K$ be a number field. Then there exists a finite set $E\subseteq K$ with $|E|\le [K\colon \Q]$ and an effective positive constant $B=B_K$ with the following property: if $\gamma\in K^c$ is an algebraic integer then there exist elements $c_1,\ldots,c_b\in E$ and roots of unity $\zeta_1,\ldots,\zeta_b\in K^c$ such that $\gamma=\sum_{i=1}^bc_i\zeta_i$, where $b\le |E|\cdot L(B\cdot \house{\gamma})$ for any Loxton function $L\colon \R\to \R$. Moreover, both $E$ and $B$ can be effectively computed in terms of a basis for $K/\Q$.
\end{theorem}

Lemma \ref{absolute_bounds} and Theorem \ref{loxton} easily imply the following corollary.

\begin{corollary}\label{conditions}
Under the hypotheses of Theorem \ref{o_theorem}, there exists a finite set $E\subseteq K$ and an effective positive constant $B$ with the following property: if $\gamma\in \phi^{-\infty}(\alpha)\cap K^c$ then there exist elements $c_1,\ldots,c_B\in E$ and roots of unity $\zeta_1,\ldots,\zeta_B\in K^c$ such that $\gamma=\sum_{i=1}^Bc_i\zeta_i$.
\end{corollary}

The third preliminary result is the following theorem of Fuchs and Zannier.

\begin{theorem}[{{\cite[Corollary]{FZ}}}]\label{fz_thm}
Let $K$ be a number field and $h\in K(x)$ be a rational function of degree $d\ge 2$ that is not $\overline{K}$-conjugate to $x^{\pm d}$ or to $\pm T_d(x)$. Then for every integer $n \geq 3$ and every non-constant $q\in K(x)$, the function $h^n(q(x))$ cannot be expressed by a ratio of polynomials having altogether less than $\frac{1}{\log 5}((n-2)\log d-\log 2016)$ terms.
\end{theorem}

Now let us explain how the proof of Theorem \ref{o_theorem} works. From now on, suppose that $\phi$ is not $\overline{K}$-conjugate to $x^{\pm d}$ or $\pm T_d(x)$, filling in the remaining cases is rather elementary. Suppose by contradiction that $\phi^{-\infty}(\alpha)\cap K^c$ is infinite. Let $M$ be a positive integer, and consider the equation
\begin{equation}\label{eq1}
\phi^M(z_1)=z_2,
\end{equation}
where $z_1,z_2\in \phi^{-\infty}(\alpha)\cap K^c$. Since the latter set is infinite, there exist infinitely many pairs $(z_1,z_2)$, with both $z_i$'s in the set, that satisfy equation \eqref{eq1}. On the other hand, by Corollary \ref{conditions} all elements of these pairs can be represented as $\sum_{i=1}^Bc_i\zeta_i$. where the $c_i$'s are chosen in a finite set $E$. By the pigeonhole principle, this means that there exist $\overline{c}_1,\ldots,\overline{c}_B\in E$ such that there are infinitely many pairs $(z_1,z_2)$ that satisfy equation \eqref{eq1} and such that both $z_1$ and $z_2$ can be expressed as $\sum_{i=1}^B\overline{c}_i\zeta_i$ for some choice of roots of unity $\zeta_1,\ldots,\zeta_B$. But this shows that the hypersurface
$$S\colon\phi^M\left(\sum_{i=1}^B\overline{c}_ix_i\right)-\sum_{i=1}^B\overline{c}_iy_i=0\subseteq \mathbb G_m^{2B}$$
has infinitely many torsion points. By the torsion coset theorem, the Zariski closure of the subset of torsion points must therefore contain a torsion coset of dimension $1$, and it is easy to see that one can choose it in such a way that it is not mapped to a constant via the map $S\to \P^1$ that sends $(\overline{x}_1,\ldots,\overline{x}_B,\overline{y}_1,\ldots,\overline{y}_B)$ to $\sum_{i=1}^B\overline{c}_i\overline{x}_i$. This means that if we parametrize such coset by writing $x_i=\zeta_it^{e_i}$ and $y_j=\zeta_jt^{e_j}$, where $\zeta_i,\zeta_j$ are roots of unity and $e_i,e_j$ are integers, not all zero, we get the following identity of rational functions:
\begin{equation}\label{eq2}
\phi^M\left(\sum_{i=1}^B\overline{c}_i\zeta_it^{e_i}\right)=\sum_{j=1}^B\overline{c}_j\zeta_jt^{e_j}.
\end{equation}
However, $M$ is arbitrary. Hence if we choose it to be large enough with respect to $B$, we obtain a representation of the rational function $\phi^M\circ q(t)$, where $q(t)=\sum_{i=1}^B\overline{c}_i\zeta_it^{e_i}$, as a rational function (namely the RHS of \eqref{eq2}) that is a ratio of polynomials with altogether at most $B+1$ nonzero terms. Since $\phi$ is not $\overline{K}$-conjugate to $x^{\pm d}$ or $\pm T_d(x)$, this contradicts Theorem \ref{fz_thm}.

Clearly, Theorem \ref{o_theorem} implies the following corollary.

\begin{corollary}
Conjecture \ref{ap_conjecture} holds true for any pair $(f,\alpha)$ where $f\in \Q[x]$ and $\alpha\in \Q$.
\end{corollary}

It is clear how powerful Theorem \ref{o_theorem} is, as it implies a far stronger result than that of Conjecture \ref{ap_conjecture}, for rational functions over the rationals. In fact, it is not unreasonable to expect that such conjecture holds true even for pairs $(\phi,\alpha)$ with $\alpha\in \Q$ and $\phi\in \Q(x)$, once one adds to the statement the obvious extra case of $(\phi,\alpha)$ being $\Q^c$-conjugate to $(x^{-d},\zeta)$. On the other hand, one cannot expect Conjecture \ref{ap_conjecture} to hold true for rational functions over any number field, as the theory of complex multiplication shows that there are infinitely many examples of Latt\`es maps (that are never conjugated to $x^{\pm d}$ or $\pm T_d(x)$) that give rise to abelian dynamical Galois groups.

\bibliographystyle{plain}
\bibliography{bibliography}

\end{document}